\documentclass[12pt]{article}
\usepackage{amsmath,amsfonts,latexsym}

\numberwithin{equation}{section}

\newcommand{\C}{{\mathbb C}}

\newcommand{\Z}{{\mathbb Z}}

\begin{document}
\begin{center}
{\Large{\bf TOROIDAL $Z$-ALGEBRAS}} \\ [3.5cm]
{\bf S. Eswara Rao} \\
{\bf School of Mathematics} \\
{\bf Tata Institute of Fundamental Research} \\
{\bf Homi Bhabha Road} \\
{\bf Mumbai - 400 005}
 {\bf India} \\ [5mm]
{\bf email: senapati@math.tifr.res.in}
\end{center}

\begin{abstract}
The toroidal Lie algebras an $N$ variable generalizations of affine Kac-Moody Lie algebras.
As in the affine Lie algebra there exists finite order automorphisms corresponding to Dynkin diagram
automorphisms. The fixed point subalgebra are called twisted toroidal Lie algebras. In this paper we 
construct faithfull representations for toroidal Lie algebras ( this includes the non-twisted case also ) useing 
methods developed by Lepowsky Wilson $[LW]$. This construction recovers the result by Eswara Rao - Moody $[EM]$ in 
the homogeneous picture and by Yuly Billig $[B1]$ in the principal picture. The proofs given in this paper are much 
shorter than above works. The results for the twisted case are completely new. 
\end{abstract}
MSC : Primary 17B67. Secondary 17B65, 17B69.\\
Keywords; Twisted toroidal Lie algebras, Vertex operators and Z-algebras.

\newpage
\noindent
{\Large{\bf Introduction}}
\vskip 5mm

In recent times a great amount of research has been done on Extended Affine Lie
Algebras (EALA), which are natural generalization of affine Kac-Moody Lie 
algebras.  See [AABGP], [AG] and references there in.    It is an accepted 
fact that the Lie algebras gains importance only when it admits  a natural
realization in other words a faithful representation.  It is an open 
problem to find a realization for an EALA.   An important class of EALA's  
are the ones obtained from the so called toroidal Lie algebras.  Toroidal
Lie algebras are $N$ variable generalization of affine Kac-Moody Lie algebras. 
 For the first time a large class of (integrable) modules are constructed for 
toroidal Lie algebras in [EM] and [MEY], the so called homogeneous picture.
In [B1] similar construction has been made for the principal picture.
These constructions are very important and have found applications in 
differential equations in the works of [B2], [ISW1] and [ISW2].  To obtain
an EALA from toroidal Lie algebra one need to  add infinite set of
derivations.  In 2006 Yuly Billig [B3] obtained a realization by 
making use of Vertex operator algebras. The next class examples of EALA's
are the one corresponding to the twisted toroidal Lie algebras.  Thus the 
purpose of this paper is to construct faithful representation of twisted
toroidal Lie algebras which  arise as fixed points of certain automorphisms 
of toroidal Lie algebras of type ADE.  The main idea
is to use the Z-algebra theory developed by Lepowsky-Wilson [LW] in the
study of Vertex operator representation for affine Kac-Moody Lie algebra.  
In the process of our construction of representations using Z-algebra
theory, we recover the results of [EM], [MEY] in the homogeneous picture
and the results of [B1] and [T] in the principle picture.  Our proofs in 
these cases are much shorter than the existing proofs.  This is the first
time we have a faithful realization for the twisted toroidal Lie algebra.

Let ${\cal G}$ be the simple finite dimensional Lie algebra over the complex
numbers.  Let $A$ be a Laurent polynomial ring in $N+1$ commuting
variables.  Consider the multiloop algebra ${\cal G} \otimes A$, its
universal central extension $\tilde{\tau}$ the toroidal  Lie algebra.
Let $\theta$ be an automorphism of ${\cal G}$ of order $m$.  Then
$\theta$ can be 
extended to an automorphism of $\tilde{\tau}$  (Section 1).  Then
the subalgebra of $\theta$ fixed points inside $\tilde{\tau}$ is called
twisted toroidal Lie algebra $\tilde{L} ({\cal G}, \theta)$.  It is the  
universal central extension of the underlining multiloop algebra (See [BK]). 

In Section 1, we define a category $\underline{C}_k $ of $\tilde{L}
({\cal G},\theta)$ modules which satisfy a factorisation property first
introduced in [BY].  The factorisation property is not satisfied for
a general class of integrable modules.  But there are enough of
integrable  modules which satisfy the factorisation property.  For
example the vertex representation defined in [EM] and  the representation
considered in [BY] satisfy factorisation property.

Next by following [LW] closely we define toroidal  $Z$-algebras (1.10) and
define a category $\underline{D}_k$-of $Z$-algebra modules.  We then
prove the important Proposition (2.6) which says that the categories
$\underline{C}_k$ and $\underline{D}_k$ are equivalent.  Thus by
constructing a $Z$- algebra module we get a module for
$\tilde{L}({\cal G},\theta)$. 

In section 3 we specialise to the homogeneous picture for the nontwisted 
case of type ADE.  We construct a module for the $Z$-toroidal Lie
algebra closely following the results of [LP].  Thereby constructing
a module for $\tilde{L} ({\cal G}, I_d) \cong \tilde{\tau}$ which  is
faithful.  This recovers the main result of [EM].  Our calculations
are certainly much shorter. 

In Section 4 we specialise to the principal picture.  This includes
the twisted and nontwisted toroidal Lie algebras.  We again
construct a module for the $Z$-toroidal Lie algebra by making use of
the corresponding results for the affine Kac-Moody Lie algebra from
[LW].  We have to consider the additional Fock space for this purpose.
Thus we get a module for our $\tilde{L} ({\cal G}, \theta)$. This result
recovers the main result of [B1] and [T].  Again our proof are much
shorter.  {\bf The twisted case is completely new.}

In the process we have given the following realization of twisted
toroidal Lie-algebra.  Let $\pi$ be a Dynkin diagram automorphism of
${\cal G}$.   Define an automorphism $\theta$ of ${\cal G}$ as in the
section 4. 
Then we prove that $\tilde{\tau} ({\cal G}, \theta) \cong \tilde{L}
({\cal G}, \pi)$. 
This is what is called the principal realization in the affine case.  The
isomorphism is given explicitly in twisted case and it is completely
new even in the affine case (Proposition 4.10). 

\section*{Section 1}

Let ${\cal G}$ be a finite dimensional  semisimple Lie-algebra over the complex
numbers $\C$.  Let $<,>$ be a non-degenerate symmetric ${\cal G}$-invariant 
bilinear form on ${\cal G}$.  We fix a non-negative integer $N$.  Let
$A=\C[ t^{\pm 1}, t_1^{\pm 1}, \cdots, t_N^{\pm 1}]$ be the ring of
Laurent polynomials in $N+1$ commuting variables.  Let
$\underline{r}=(r_1, \cdots, r_N) \in \Z^N$.  Let $t^{\underline{r}} =
t_1^{r_1} t_2^{r_2} \cdots t_N^{r_N}$. Fix a positive integer $m$.  
Let $\Omega_A$ be a free
$A$-module of rank $N+1$  with basis $\{k_0, \cdots, k_N\}$.  Let
$d_A$ be the subspace of $\Omega_A$ spanned by elements of the form
$\frac{1}{m} r_0 t^{r_0} t^{\underline{r}} k_0+ \cdots + r_N t^{r_0}
t^{\underline{r}} k_N$.  Let $x (r_0, \underline{r}) = x \otimes
t^{r_0} t^{\underline{r}} \in {\cal G} \otimes A$.  Then the toroidal
Lie-algebra $\tau = {\cal G} \otimes A \oplus \Omega_A /d_A$ is defined by
the following bracket.

$$
\begin{array}{lll}
[x (r_0,\underline{r}), y (s_0, \underline{s})] &=& \\
\left[ x,y \right] (r_0+s_0, \underline{r}+ \underline{s})+ \frac{<
x,y>}{m} r_0 
t^{r_0+s_0} t^{\underline{r}+ \underline{s}}k_0 
&+& <x,y> \displaystyle{\sum_{i=1}^{N}} r_i t^{r_0+s_0}
t^{\underline{r}+ \underline{s}} k_i
\end{array}
$$
 for $x,y \in {\cal G},
\underline{r}, \underline{s} \in \Z^N, r_0, s_0 \in \Z$.\\
 $\Omega_A/d_A$ is central.  

It is known that $\tau$ is the universal central extension of ${\cal G}
\otimes A$. (See [K], [MEY]).  (First note that the toroidal Lie
algebra defined by $m=1$ is isomorphic to the above).   Let $\underline{h}$ be a
Cartan subalgebra of ${\cal G}$.  Let $\theta$ be an automorphism of
${\cal G}$ such that 
$\theta (\underline{h}) =\underline{h}$ and of order $m$.  
 Let
$\Z_m=\Z/_{m  \Z}$ 
be the cyclic group of order $m$.  Let $w$ be a primitive $m$ th root of
unity. 

\paragraph*{(1.2)} Let ${\cal G}_i = \{x \in {\cal G} \mid
\theta x =w^ix \}$ for $i \in 
\Z$.  Then ${\cal G} = \displaystyle{\oplus_{i \in \Z_m}} {\cal G}_i$.
Note that $< 
{\cal G}_i, {\cal G}_j > =0$ unless $i+j \equiv 0 (m)$.  For $x
\in {\cal G}$ write 
$x=\displaystyle{\sum_{i \in \Z_m}} x_i$ where $\theta x_i=w^i x_i$.
Define $x_i=x_{\overline{i}}$ for $ i \in \Z$ and $\overline{i}
\in \Z_m$.  

Extend the automorphism $\theta$ to $\tau$ by $\theta (x
(r_0,\underline{r}))= {w}^{-r_o} \theta (x) (r_0, \underline{r})$
and $\theta (t_0^{r_o} t^{\underline{r}} k_i) = {w}^{-r_0}
t^{r_0} t^{\underline{r}} k_i, 0 \leq i \leq N$.  Let $\tilde{\tau}
= \tau \oplus D$ where $D$ is spanned by derivations $\{d_0, \cdots,
d_N\}$ with bracket $[d_i, x(r_0, \underline{r})] =r_i x (r_0,
\underline{r})$, for $0 \leq i \leq N$, 
$[d_i, t^{ro} t^{\underline{r}}k_{j}]= r_i t^{ro} t^{\underline{r}}k_{j}$ and $[d_i, d_j]=0.$
Extend the automorphism $\theta$ to $\tilde{\tau}$ by  $\theta (d_i) =d_i$.  Let $(\Omega_A / {d_A})_0$ be
the linear span of $t^{r_0} t^{\underline{r}} k_i$ where $r_0
\equiv 0 (m)$.  Consider the $\theta$ fixed points of $\tilde{\tau }$ say
$\tilde{L} ({\cal G}, \theta)$. 

\paragraph*{(1.3)}.  Let $L ({\cal G}, \theta) =
\displaystyle{\oplus_{\stackrel{i \in \Z}{\underline{r} \in \Z^N}}}
{\cal G}_i (i, \underline{r})$ and $\overline{L} ({\cal G}, \theta) =
L ({\cal G}, \theta) 
\oplus (\Omega_A / {d_A})_0$.  Then clearly $\tilde{L} ({\cal G}, \theta) =
\overline{L} ({\cal G}, \theta) \oplus D$.

Since $<,>$ is non-degenerate and ${\cal G}$-invariant, its restriction to
$\underline{h}$ is also non-degenerate.  We identify $\underline{h}$
and $\underline{h}^*$ via this form.  Let $\Phi$
be the root system of ${\cal G}$.  For $\beta \in \Phi$, choose the
corresponding non-zero root vectors $x_{\beta}$ such that
$[x_{\beta}, x_{- \beta}] =< x_{\beta},
{x}_{- \beta} > \beta$.  Let $\epsilon (\beta, \gamma)$ be a
non-zero number such that $[x_{\beta}, x_{\gamma}]= \epsilon (\beta,
\gamma) x_{\beta+\gamma}$.  Clearly the set of roots $\Phi $ is $\theta$-
stable.  Then define $\eta (p, \beta)$ a non-zero scalar such that

\paragraph*{(1.4)} $\theta^p x_{\beta} = \eta (p, \beta) x_{\theta^p \beta}$.

For any vector space $V$ and for indeterminates
$\zeta_1, \cdots, \zeta_\ell$, denote $V \{\zeta_1, \cdots,
\zeta_{\ell}\}$ the space of formal Laurent series.  Further $V
[\zeta_1^{\pm 1}, \cdots, \zeta_n^{\pm 1}]$ denote finite formal Laurent
series.  We recall the following Proposition from [LW].  Define
$\delta (\zeta) = \displaystyle{\sum_{i \in \Z}} \zeta^i \in
\C\{\zeta\}$.

\paragraph*{Proposition ~(1.5) (a)} (Proposition (2.2) of [LW]).  Let $f
(\zeta) \in V [\zeta, \zeta^{-1}]$. Then
$$f (\zeta) \delta (\zeta^m) = m^{-1}\displaystyle{\sum_{p \in \Z_m}} f (w^p)
\delta (w^{-p} \zeta).$$
For $x \in {\cal G}, \underline{r} \in \Z^N$
let $$x (\underline{r}, \zeta) = \sum_i x_i \otimes t^i t^{\underline{r}}
\zeta^i, \ k_i (\underline{r}, \zeta^m) = \sum_{p \in \Z} t^{mp} t^r k_i
\zeta^{mp}.$$ 
For any infinite series $f (\zeta) = \sum b_i \zeta^i$, let $Df (\zeta)= \sum i b_i \zeta^i$. 

\paragraph*{Proposition (1.5)}(b) The following relations hold for $x
(\underline{r}, \zeta)$ and $k_i (\underline{r}, \zeta^m)$.  In fact
they define a Lie-algebra $\tilde{L} ({\cal G}, \theta)$.  For
$\beta_1, \beta_2 
\in \Phi, \underline{r}, \underline{s} \in \Z^N$.
$$[x_{\beta_1} (\underline{r}, \zeta_1), x_{\beta_2} (\underline{s},
\zeta_2)] \leqno{(1)}$$
$$= \frac{1}{m} \displaystyle{\sum_{\theta^p \beta_1+ \beta_2 \in 
\Phi}} \eta (p, \beta_1) \epsilon (\theta^p \beta_1, \beta_2) x_{\theta^p
\beta_1+ \beta_2} 
(\underline{r}+\underline{s}, \zeta_2) \delta (w^{-p} \zeta_1/\zeta_2)$$ 
$$- \frac{1}{m} < x_{\beta_2}, {x}_{- \beta_2} >
\displaystyle{\sum_{\theta^p \beta_1+ \beta_2 =0}} \eta
(p, \beta_1) \beta_2 (\underline{r}+ \underline{s}, \zeta_2) \delta
(w^{-p} \zeta_1/\zeta_2)$$

$$+ \frac{1}{m} < x_{\beta_2}, x_{- \beta_2} >
\displaystyle{\sum_{\theta^p \beta_1+ \beta_2 =0}} \eta
(p, \beta_1) (\displaystyle{\sum_{i=1}^{N}} r_i k_i (\underline{r}+
\underline{s}, \zeta_2^m) \delta (w^{-p} \zeta_1/\zeta_2)$$
$$+\frac{k_0(\underline{r}+ \underline{s}, \zeta^m_2)}{m} D \delta
({w}^{-p} \zeta_1/\zeta_2))$$
$$[\beta_1 (\underline{r}, \zeta_1), \beta_2 (\underline{s},
\zeta_2)] \leqno{(2)}$$
$$= \frac{1}{m} \displaystyle{\sum_{p \in \Z_m}} < \theta^p \beta_1,
\beta_2> ( \displaystyle{\sum_{i=1}^{N}} r_i k_i
(\underline{r}+\underline{s}, \zeta_2^m) \delta ({w}^{-p}
\zeta_1/ \zeta_2)$$
$$+ \frac{k_0 (\underline{r}+ \underline{s}, \zeta^m_2)}{m} D \delta
(w^{-p} \zeta_1/ \zeta_2)).$$
$$[\beta_1 (\underline{r}, \zeta_1), x_{\beta_2} (\underline{s},
\zeta_2)]= \leqno{(3)}$$
$$\frac{1}{m} \displaystyle{\sum_{p \in \Z_m}} < \theta^p \beta_1,
\beta_2> x_{\beta_2} (\underline{r}+ \underline{s}, \zeta_2) \delta
(w^{-p} \zeta_1/\zeta_2)$$
$$\frac{1}{m} D k_0 (\underline{r}, \zeta^m)+
\displaystyle{\sum_{i=1}^{N}} r_i k_i (\underline{r}, \zeta^m) =0
\leqno{(4)}$$ 
$$[d_i, k_j (\underline{r}, \zeta^m)] = r_i k_j (\underline{r},
\zeta^m) \ {\rm for} \ 1 \leq i \leq N \ {\rm and} \ \ 0 \leq j \leq
N \leqno{(5)}$$ 
$$[d_0, k_j (\underline{r}, \zeta^m)]=D k_j (\underline{r}, \zeta^m),
\ 0 \leq j \leq N
\leqno{(6)}$$ 
$$x_{\beta} (\underline{r}, w^p \zeta) = \eta(p, \beta) x_{\theta^p
\beta} (\underline{r}, \zeta) \leqno{(7)}$$
\[k_i (\underline{r}, \zeta^m) \mbox{ is \ central  in } \overline{L}({\mathcal G},\theta)  \mbox{ for  } 
0 \leq i \leq N \leqno{(8)}\]
$$[d_i, x_\beta (\underline{r}, \zeta)] = r_i x_\beta (\underline{r}, \zeta), \ 1 \leq i \leq N \leqno{(9)}$$
$$[d_0, x_\beta (\underline{r}, \zeta)] = D x_\beta (\underline{r}, \zeta), \ [d_i, d_j]=0 \leqno{(10)}$$

\paragraph*{Proof}  (4)  follows from the definition of
$(\Omega_A/d_A)_0$. (5), (6), (7) and (8) are easy to see.  
First consider the following:
$$[x (\underline{r}, \zeta_1), y (\underline{s}, \zeta_2)]=F+G_1+G_2.$$
Where 
$$F= \displaystyle{\sum_{i,j}}  [x_i, y_j] t^{i+j}
t^{\underline{r}+ \underline{s}}  \zeta_1^i \zeta_2^j,$$
$G_1=\frac{1}{m} \ \displaystyle{\sum_{i,j}} \ i < x_i,
y_j> t^{i+j} t^{\underline{r}+ \underline{s}} k_0   \zeta_1^i
\zeta_2^j$,
$$G_2 = \displaystyle{\sum_{\ell=1}^{N}} \displaystyle{\sum_{i,j}} r_{\ell}
<x_i, y_j > t^{i+j} t^{\underline{r}+\underline{s}} k_\ell \zeta_1^i
\zeta_2^{j}.$$
From the proof of Theorem (2.3) of [LW] it follows that
$$F= \frac{1}{m} \displaystyle{\sum_{p \in \Z_m}} [\theta^p x, y]
(\underline{r}+ \underline{s}, \zeta_2) \delta (w ^{-p} \zeta_1 /
\zeta_2).$$
Now consider $ [\theta^p x_{\beta_1}, x_{\beta_2}]= \eta (p, \beta_1) [x_{\theta^p \beta_1}, x_{\beta_2}]$.\\
$=\begin{cases} \eta (p, \beta_1) \epsilon (\theta^p \beta_1, \beta_2) x_{\theta^
p \beta_1+\beta_2} \ {\rm if} \ \theta^p \beta_1+ \beta_2 \in \Phi\\
-\eta (p, \beta_1)<x_{\beta_2}, x_{- \beta_2} > \beta_2 \ {\rm if} \ \theta^p \beta_1+ \beta_2=0 \\
0 \ \mbox{if}  \  \theta^{\beta} \beta_1+ \beta_2 \neq 0 \  \ \mbox{and \ not \ a \ root.}\end{cases}$\\[2mm]   

Thus for $x=x_{\beta_1}$ and $y= x_{\beta_2}, F$ equals to the first and 
second term of right hand side in (1).  For $G_1$ and $G_2$, first note that
$<x_i, y_j> =0$ if $i+j \not\equiv 0 (m)$ by (1.2).  Thus
$$mG_1= \displaystyle{\sum_{i,p}} i<x_i, y_{-i+mp}> t^{mp}
t^{\underline{r}+\underline{s}} k_0 (\zeta_1/ \zeta_2)^i \zeta_2^{mp}.$$
$$=\displaystyle{\sum_{i,p}} i < x_i, y_{-i} > t^{mp}
t^{\underline{r}+\underline{s}} k_0 (\zeta_1/ \zeta_2)^i \zeta_2^{mp}$$
$$=\displaystyle{\sum_{i,p}} i <x_i, y> t^{mp}
t^{\underline{r}+\underline{s}} k_0 (\zeta_1/ \zeta_2)^i \zeta_2^{mp}$$
$$=D (\displaystyle{\sum_{i}} <x_i, y> (\zeta_1/ \zeta_2)^i) k_0
(\underline{r}+\underline{s}, \zeta_2^m)$$
$$=D(\displaystyle{\sum_{i=0}^{m-1}} <x_i, y> (\zeta_1/\zeta_2)^i \delta
((\zeta_1 / \zeta_2)^m).k_0 (\underline{r}+\underline{s}, \zeta_2^m)$$
$$= \frac{1}{m} \displaystyle{\sum_{p}} < \theta^p x,y>D \delta(w^{-p}
\zeta_1/\zeta_2) k_0 (\underline{r}+\underline{s}, \zeta_2^m). \eqno{1.5(1)} $$
Similarly
$$G_2= m^{-1}\displaystyle{\sum_{\ell=1}^{N}} r_{\ell} < \theta^p x,y> \delta 
(w^{-p} \zeta_1/ \zeta_2) k_{\ell} (\underline{r}+ \underline{s},
\zeta_2^m ) ~ \eqno{  1.5(2)}$$

Again for $x=x_{\beta_1}$ and $y=x_{\beta_2}, < \theta^p x, y> = \eta (p,
\beta_1) < x_{\theta^p \beta_1}, x_{\beta_2}>=0 \ {\rm if} \ \theta^p
\beta_1+\beta_2 \neq 0$ \\
$=\eta(p, \beta_1) <x_{-\beta_2},x_{\beta_2}> \ {\rm if} \ \theta^p
\beta_1+\beta_2=0$. \\
This completes the proof (1).  To see (2), take $x=\beta_1$ and
$y=\beta_2$.  Then $F=0$.  From 1.5 (1) and 1.5(2), (2) will follow.
To see (3) take $x=\beta_1$ and $y= x_{\beta_2}$ and note that $G_1= 0$
and $G_2=0$ and $F$ is equal to right hand side of 3.  This completes the proof of the 
Proposition (1.5) (b).

Now we define a category $\underline{C}_k $ of $\tilde{L} ({\cal G},
\theta)$- modules.  

\paragraph*{Definition (1.6)} A $\tilde{L} ({\cal G}, \theta)$-module $V$ is 
in $\underline{C}_k$ if
\begin{enumerate}
\item[(1)] $k_o$ acts as $k$.
\item[(2)] $V= \displaystyle{\oplus_{z \in \C}} V_z, V_z = \{ v \in V
\mid d_0 v=zv \}$.  Assume for any $z$ there exists $\ell_0$ such that 
$ V_{z+ \ell}
=0$ for $ \ell > \ell_0$ 
\item[(3)] $\frac{1}{k} x ( \underline{r}, \zeta) k_0 (\underline{s},
\zeta^m) =x (\underline{r}+ \underline{s}, \zeta) $ \\
$\frac{1}{k} k_i (\underline{r}, \zeta^m) k_o (\underline{s},
\zeta^m) = k_i (\underline{r}+ \underline{s}, \zeta^m) \ {\rm for} \
0 \leq i \leq N$.
\end{enumerate}

\paragraph*{Remark (1.7)}  Condition (3) is not satisfied for most
of the modules.  But there are enough of them which are sufficient for
a realization of $L ({\cal G}, \theta)$.  For examples vertex operator
representation of [EM] satisfy the condition (3) as well as the representations
considered in [BY].

\paragraph*{(1.8)} Consider the Lie subalgebra of
$\tilde{L} ({\cal G}, \theta), \tilde{\underline{h}} =
\displaystyle{\oplus_{i \in \Z}} h_i \otimes t^i \oplus \C k_0 \oplus
\C d_0 (h \in
\underline{h})$. The bracket is given by 
$$[h_i \otimes t^i, h_j \otimes t^j] = \frac{i k_0}{m}  <h_i, h_j >
\delta_{i+j,0}.$$ 
Clearly $\tilde{\underline{h}}$ is
$\Z$-graded.  Let $M (k)$ be a Verma module of level $k$ for
$\tilde{\underline{h}}$.  Then it is a standard fact that $M(k)$ is
irreducible whenever $k$ is non-zero.

\paragraph*{Proposition (1.9)}  Any module $V $ in $\underline{C}_k (k
\neq 0)$ has the following decomposition as
$\tilde{\underline{h}}$-modules. 
$V \cong M (k) \otimes \Omega_V$ where 
$$\Omega_V = \{v \in V \mid h_i \otimes t^i v=0 \  {\rm for} \ i >0 \}$$
see Proposition 5.4 of [LW].

We now define toroidal $Z$ algebras. Notation as earlier.  For
$\alpha \in \Phi, \underline{r} \in \Z^N$ let $Z (\alpha,
\underline{r}, \zeta)$ be a series in $\zeta $ with integral powers.
For $\underline{r} \in \Z^N$ let $k_i (\underline{r}, \zeta^m)$ be a
series in $\zeta^m$.  The toroidal $Z_k$-algebra or simply
$Z_k$-algebra is an algebra generated  by the components of $Z (\alpha,
\underline{r}, \zeta), k_i (\underline{r}, \zeta^m), \underline{h}_0$ and $d_0, d_1, \cdots, d_N$
by the following relations. 

\paragraph*{(1.10) ~~ Relations} $\alpha, \beta, \beta_1, \beta_2 \in
\Phi, \underline{r}, \underline{s} \in \Z^N$.

\begin{enumerate}
\item[(1)] $\frac{1}{k} Z (\alpha, \underline{r}, \zeta) k_0
(\underline{s}, \zeta^m) = Z (\alpha, \underline{r}+ \underline{s},
\zeta)$ 
\item[(2)] $\frac{1}{k} k_o( \underline{r}, \zeta^m) \ k_i
(\underline{s}, \zeta^m) = k_i ( \underline{s}+ \underline{r},
\zeta^m), \  0 \leq i \leq N$ .
\item[(3)] $\displaystyle{\sum_{i=1}^{N}} r_i k_i (\underline{r},
\zeta^m)+ \frac{1}{m} D k_0 (\underline{r}, \zeta^m) =0$.
\item[(4)] $[d_0, Z (\beta, \underline{r}, \zeta)]=D Z (\beta,
\underline{r}, \zeta)$ 
\item[(5)] $[d_0, k_i (\underline{r}, \zeta^m)] = D k_i
(\underline{r}, \zeta^m), \ 0 \leq i \leq N$ 
\item[(6)] $[d_i, k_j (\underline{r}, \zeta^m)]= r_i k_j
(\underline{r}, \zeta^m), \ 1 \leq i \leq N, 0 \leq j \leq N$
\item[(7)] $ \displaystyle{\prod_{p \in \Z_m}} (1-w^{-p}
\zeta_1/\zeta_2)^{< \theta^p \beta_1, \beta_2>/k} Z (\beta_1, \underline{r}, \zeta_1) Z (\beta_2,
\underline{s}, \zeta_2)$ \\
$-\displaystyle{\prod_{p \in \Z_m}} (1-w^{-p}
\zeta_2/\zeta_1)^{<\theta^p \beta_2, \beta_1>/k} Z 
(\beta_2, \underline{s}, \zeta_2) Z(\beta_1, \underline{r}, \zeta_1)$ \\
$=\frac{1}{m} \displaystyle{\sum_{\theta^p \beta_1+ \beta_2 \in
\Phi}} \eta (p, \beta_1) \epsilon  (\theta^p \beta_1, \beta_2) Z
(\theta^p \beta_1+\beta_2, \underline{r}+\underline{s}, \zeta_2)
\delta (w^{-p} \zeta_1/ \zeta_2)$ \\
$- \frac{1}{mk} < x_{\beta_2}, {x}_{-\beta_2}>
\displaystyle{\sum_{\theta^p \beta_1+ \beta_2=0}} \ \eta (p,
\beta_1) (\beta_2)_0 k_0 (\underline{r}+ \underline{s},
\zeta_2^m) \delta (w^{-p} \zeta_1/\zeta_2)$ \\
$+ \frac{1}{m} <x_{\beta_2}, {x}_{-\beta_2} > \
\displaystyle{\sum_{\theta^p \beta_1+ \beta_2=0}} \ \eta (p, \beta_1)
(\displaystyle{\sum_{i=1}^{N}} \ r_i k_i (\underline{r}+\underline{s}, \zeta_2^m)
\delta (w^{-p} \zeta_1/ \zeta_2)+\frac{1}{m} k_0
(\underline{r}+\underline{s}, \zeta_2^m) D
\delta (w^{-p} \zeta_1/ \zeta_2))$.
\item[(8)] $[\alpha, Z (\beta, \underline{r}, \zeta)] = <\alpha,
\beta > Z (\beta, \underline{r}, \zeta), \ \alpha \in \underline{h}_0$ 
\item[(9)] $Z (\beta, \underline{r}, w^{p} \zeta) = \eta(p, \beta) Z
(\theta^p \beta, \underline{r}, \zeta)$
\item[(10)] $k_i (\underline{r}, \zeta^m)$ commutes with $Z(\alpha, \underline{r}, \zeta)$ and 
$\underline{h}_0$,  \ $0 \leq i \leq N$. 
\end{enumerate}

As it is we do not know whether a $Z_k$ algebra is non-zero or not but
certainly it is well defined. 

\paragraph*{(1.11) Definition} $A ~ Z_k$ module $V$ is said to be in
the category $\underline{D}_k$ if \\
(1) $ k_0$ acts by scalar $k$. \\
(2) $V = \displaystyle{\oplus_{z \in \C}} V_z, V_z = \{v
\in V \mid d_0 
v=z v \}$. Assume for any given $z$ there exists $ \ell_0$ such that 
 $V_{z+ \ell} = 0$ for $\ell > \ell_0$.

\section*{Section 2} 

In this section we establish equivalence between the categories
$\underline{C}_k$ and $\underline{D}_k$.  The proof are very similar to
[LW].  In fact most of the results go through.

Let $V \in \underline{C}_k$.  Define for $\beta \in \Phi$.
$$E^{\pm} (\beta, \zeta) =  {\rm exp} (\pm m \sum_{j>0} \beta_{\pm j} \otimes
t^{\pm j} \zeta^{\pm j}/jk)$$
and $Z (\beta, \underline{r}, \zeta) = E^- (\beta, \zeta) x_{\beta}
(\underline{r}, \zeta) E^+ (\beta, \zeta)$.  We first prove that
these $Z$-operators satisfy relations in (1.10).

We first recall the following from section 3 of [LW].  We are taking
$\underline{a} = \underline{h}$ and $\underline{m} =0 $ (in decomposition 
(3.5)) in [LW].

\paragraph*{Proposition (2.1)} $\alpha, \beta, \gamma \in \Phi$.
Let
\begin{enumerate}

\item[(1)]~~~~~  $i>0$
\item[(a)] $ [\alpha_i \otimes t^i, E^+ (\beta, \zeta)] =0$
\\
\item[(b)] $[\alpha_{-i} \otimes t^{-i}, E^- (\beta, \zeta)] =0$ \\
\item[(c)] $[ \alpha_i \otimes t^i, E^- (\beta, \zeta)]= - < \alpha_i,
\beta > \zeta^{-i} E^- (\beta, \zeta)$. \\
\item[(d)] $[ \alpha_{-i} \otimes t^{-i}, E^+ (\beta, \zeta) ]
=  -< \alpha_{-i}, \beta > \zeta^i E^+ (\beta, \zeta)$ \\
\item[(2) (a)] $ ~~~ E^{\pm} (\beta+ \gamma, \zeta) =  E^{\pm} (\beta,
\zeta) E^{\pm}(\gamma, \zeta) $\\
\item[(b)] $E^{\pm} (\theta^p \beta, \zeta)  =  E^{\pm} (\beta, w^p
\zeta)$ \\
\item[(c)] $D E^{\pm} (\beta, \zeta) = \frac{m}{k} \beta (\zeta)^{\pm}
E^{\pm} (\beta, \zeta)$ \\
${\rm where} \ \beta (\zeta)^{\pm}  = \sum_{i>0} \beta_{\pm i} \otimes
t^{\pm i} \zeta^{\pm i}$ \\
\item[(3)] (a) ~~$  x_{\beta} (\underline{r}, \zeta) =E^- (- \beta, \zeta)
Z (\beta, \underline{r}, \zeta) E^+ (- \beta, \zeta)$ \\
\item[(b)] $Z (\beta, \underline{r}, w^p \zeta) = \eta (p, \beta) Z
(\theta^p \beta,
\underline{r}, \zeta)$ \\
\item[(4)] ~~~~~ $E^+ (\beta, \zeta_1) E^- (\gamma, \zeta_2) = E^-
(\gamma,\zeta_2) E^+ (\beta, \zeta_1)$.
 $\displaystyle{\prod_{p \in \Z_m}} (1- w^{-p} \zeta_1 / \zeta_2)^{ \frac{< \theta^p \beta, \gamma>}{k}}$ \\
\item[(5)] ~~~ $E^+ (\beta, \zeta_1) x_{\gamma} (\underline{r}, \zeta_2)
= x_{\gamma} (\underline{r}, \zeta_2) E^+ (\beta, \zeta_1)$
$\displaystyle{\prod_{p \in \Z_m}} (1- w^{-p} \zeta_1/\zeta_2)^{\frac{-< \theta^p \beta, \gamma >}{k}}$ \\
\item[(6)] $ x_{\beta} (\underline{r}, \zeta_1) E^- (\gamma, \zeta_2) =E^-
(\gamma, \zeta_2) x_{\beta} (\underline{r}, \zeta_1) 
\displaystyle{\prod_{p \in \Z_m}} (1- w^{p} \zeta_1/ \zeta_2)^{\frac{-<
\theta^p \beta, \gamma >}{k}}$.

\end{enumerate}
\paragraph*{Proof} (1) a and b follows from the definition.  To see 1(c)
consider
$$[\alpha_i \otimes t^i, -\displaystyle{\sum_{j<0}} \beta_j \otimes t^j /
j_k \zeta^j]$$
$$= \frac{1}{m} < \alpha_i, \beta_{-i}> \zeta^{-i}.$$
Now
$$[\alpha_i \otimes t^i, \frac{(-\sum \beta_j \otimes t^j/j_k \zeta^j)^{\ell}}{\ell !}]$$
$$=-\frac{< \alpha_i, \beta_{-i}>}{m(\ell -1)!} \zeta^{-i} \Big(\sum \frac{\beta_j 
\otimes t^j}{j_k} \zeta^j)^{\ell-1}\Big).$$
First note that $< \alpha_i, \beta_{-i} > =<\alpha_i, \beta>$ by (1.2).  Now 1(c)
follows from the definition. \\
(1)(d) follows from similar  argument.\\
(2) and (3) follows from definition. See also Proposition 3.2 and 3.3 of [LW]. \\
(4) follows from Proposition 3.4 of [LW]. \\
(5) and (6) follows from Proposition 3.5 and 3.6 of [LW]. \\
\paragraph*{Corollary (2.2)}  Let $\underline{\tilde{h}}^{1} =
\displaystyle{\oplus_{i \neq 0}} \underline{h}_i \otimes t^i$.  Then
as operators the following hold.
(1) ~~ $[ \tilde{\underline{h}}^1, Z (\alpha, \underline{r}, \zeta)] = 0$ \\
(2) ~~ $ [a_0, Z (\alpha, \underline{r}, \zeta)]=  < a_0,
\alpha > Z (\alpha, \underline{r}, \zeta),  a \in \underline{h}$.

\paragraph*{Proof}  (1) Follows from above. (2) is easy to see.

\paragraph*{Proposition (2.3)} (Proposition (3.9) of [LW]).

Let $W$ be a vector space and let $f (\zeta_1, \zeta_2) = \sum w_{ij}
\zeta_1^i \zeta_2^j$ where  each $w_{ij} \in W$ and suppose for some
$n \in \Z$ either $w_{ij}=0$ where one of $i$ or $j> n$ or $w_{ij} =0$
whenever $i$ or $j <n$.

Set
$$D_i f (\zeta_1, \zeta_2)= \zeta_i \frac{df}{d \zeta_i} (\zeta_1,
\zeta_2)$$
Then for $a \neq 0$
$$
\begin{array}{lll}
(1) ~~~ \delta (a \zeta_1/ \zeta_2) f (\zeta_1, \zeta_2) &=& \delta
(a \zeta_1 / \zeta_2) f (\zeta_1, a \zeta_1) \\
&=& \delta ( a \zeta_1/\zeta_2) f (a^{-1} \zeta_2, \zeta_2) \\
(2) ~~ D \delta (a \zeta_1/\zeta_2) f (\zeta_1, \zeta_2) &=& (D
\delta) (a \zeta_1 / \zeta_2) f (\zeta_1, a \zeta_1)
+ \delta (a \zeta_1/\zeta_2) D_2 f (\zeta_1, a \zeta_1) \\
&=& D \delta (a \zeta_1 / \zeta_2) f (a^{-1}  \zeta_2, \zeta_2) -
\delta (a \zeta_1/\zeta_2) D_1 f (a^{-1} \zeta_2, \zeta_2) 
\end{array}
$$

\paragraph*{Proposition  (2.4) } 
$$
\begin{array}{lll}
&& \displaystyle{\prod_{p \in \Z_m}} (1-w^{-p}
\zeta_1/\zeta_2)^{< \theta^p 
\beta_1, \beta_2 > / k} Z (\beta_1, \underline{r}, \zeta_1) Z
(\beta_2, \underline{s}, \zeta_2) \\
&&-\displaystyle{\prod_{p \in \Z_m}} (1-w^{-p}
\zeta_2/\zeta_1)^{< \theta^p 
\beta_2, \beta_1 > / k} Z (\beta_2, \underline{s}, \zeta_2) Z
(\beta_1, \underline{r}, \zeta_1) \\
&=& E^- (\beta_1, \zeta_1) E^- (\beta_2, \zeta_2) [x_{\beta_1}
(\underline{r}, \zeta_1), x_{\beta_2} (\underline{s}, \zeta_2)] \\
&&E^+ (\beta_1, \zeta_1) E^+ (\beta_2, \zeta_2)
\end{array}
$$
\paragraph*{Proof}  Consider
$$Z(\beta_1, \underline{r}, \zeta_1) Z (\beta_2, \underline{s}, \zeta_2)$$
$$=E^- (\beta_1, \underline{r}, \zeta_1) x_{\beta_1} (\underline{r}, \zeta_1) 
E^+ (\beta_1, \underline{r}, \zeta_1).$$
$$E^- (\beta_2, \underline{s}, \zeta_2) x_{\beta_2} (\underline{s}, \zeta_2)
E^+ (\beta_2, \underline{s}, \zeta_2)$$
$$= \displaystyle{\prod_{p \in \Z_m}} (1- w^{-p} \zeta_1 / \zeta_2)^{\frac{< 
\theta^{p} \beta_1, \beta_2>}{k}}.$$
$$E^- (\beta_1, \underline{r}, \zeta_1) x_{\beta_1} (\underline{r}, \zeta_1)
E^- (\beta_2, \underline{s}, \zeta_2).$$
$$E^+ (\beta_1, \underline{r}, \zeta_1) x_{\beta_2} (\underline{s}, \zeta_2)
E^+ (\beta_2, \underline{s}, \zeta_2)$$
from 4 of Proposition (2.1).
$$= \displaystyle{\prod_{p \in \Z_m}} (1-w^{-p} \zeta_1/ \zeta_2)^{\frac{-< \theta^{p} \beta_1, \beta_2>}{k}}.$$
$$E^- (\beta_1, \underline{r}, \zeta_1) E^- (\beta_2, \underline{s}, \zeta_2)
x_{\beta_1} (\underline{r}, \zeta_1).$$
$$x_{\beta_2} (\underline{s}, \zeta_2) E^+ (\beta_1, \underline{r}, \zeta_1)
E^+ (\beta_2, \underline{s}, \zeta_2)$$
\Big(by 5 of Proposition (2.1)\Big).

Multiplying both sides by the inverse of the first factor on the right, and 
subtracting the expression obtained by interchanging the roles of the
subscripts 1 and 2 we have Proposition  (2.4).

\paragraph*{Proposition (2.5)}  For these $Z$ operators the relation at
(1.10) hold.

\paragraph*{Proof}  (2) to (6), holds from definition of $d_i$.
Since $V \in \underline{C}_k$ we have
$$\frac{1}{k} x_{\beta} (\underline{r}, \zeta) k_0 (\underline{s},
\zeta^m) = x_{\beta} (\underline{r}+ \underline{s}, \zeta).$$
Thus (1) holds from  definition of $Z$ operator. (9) holds from
Proposition 2.1(3). (6), (10)  and (8) are easy to see.  We only need to prove (7).
The right hand side of the Proposition (2.4) and by using Proposition 1.5 (b) (1)
 is equal to
$$E_1+ E_2 + E_3+E_4$$
$$
\begin{array}{lll}
{\rm where} \ E_1&=&E^- (\beta_1, \zeta_1) E^- (\beta_2, \zeta_2). \\
&&\frac{1}{m} \displaystyle{\sum_{\theta^p \beta_1+ \beta_2 \in
\Phi}} \ \eta (p, \beta_1) \epsilon ( \theta^p \beta_1, \beta_2)
x_{\theta^p \beta_1+ 
\beta_2}  (\underline{r}+ \underline{s}, \zeta_2) \delta (w^{-p}
\zeta_1/\zeta_2) \\
&&E^+ (\beta_1, \zeta_1) E^+ (\beta_2, \zeta_2)\\
E_2 &=&E^- (\beta_1, \zeta_1) E^- (\beta_2, \zeta_2).\\
&& (\frac{-1}{m}) <x_{\beta_2}, x_{- \beta_2} >
\displaystyle{\sum_{\theta^p \beta_1+ \beta_2 =0}} \eta (p, \beta_1)
\beta_2 (\underline{r}+ \underline{s}, \zeta_2) \delta (w^{- p} \zeta_1 /
\zeta_2). \\
&&E^+ (\beta_1, \zeta_1) E^+ (\beta_2, \zeta_2)\\
E_3&=& E^- (\beta_1, \zeta_1) E^- (\beta_2, \zeta_2). \\
&&\frac{1}{m} < x_{\beta_2}, x_{-\beta_2} >
\displaystyle{\sum_{\theta^p \beta_1+\beta_2=0}} \eta (p, \beta_1) \displaystyle{\sum_{i=1}^{N}} r_i k_i
(\underline{r}+\underline{s} , \zeta_2^m) \delta (w^{-p}
\zeta_1/\zeta_2) \\
&&E^+ (\beta_1, \zeta_1) E^+ (\beta_2, \zeta_2) \\
E_4 &=& E^- (\beta_1, \zeta_1) E^- (\beta_2, \zeta_2) \\
&&\frac{1}{m^2} < x_{\beta_2}, x_{-\beta_2}>
\displaystyle{\sum_{\theta^p \beta_1+ 
\beta_2=0}} \eta (p, \beta_1) k_0 (\underline{r}+ \underline{s},
\zeta_2^m) D \delta (w^{-p} \zeta_1/\zeta_2) \\
&&E^+ (\beta_1, \zeta_1) E^+  (\beta_2, \zeta_2)
\end{array}
$$
Now

$$E_1= \frac{1}{m} \displaystyle{\sum_{\theta^p \beta_1+ \beta_2 \in \Phi}} 
\Big(\eta (p, \beta_1)E^- (\theta^p \beta_1+ \beta_2, \zeta_2).$$
$$\epsilon(\theta^p \beta_1+\beta_2)(x_{\theta^p \beta_1+\beta_2} (\underline{r}+ \underline{s}, \zeta_2) \delta
 (w^{-p} \zeta_1/\zeta_2).$$
$$E^+ (\theta^p \beta_1+\beta_2, \zeta_2).$$
$$E_2= - \frac{1}{m} <x_{\beta_2}, x_{- \beta_2} > \displaystyle{\sum_{
\theta^p \beta_1+ \beta_2=0}} \eta (p, \beta_1) \beta_2 (\underline{r}+
\underline{s}, \zeta_2) \delta (w^{-p} \zeta_1 /\zeta_2)$$
$$E_3 = \frac{1}{m} <x_{\beta_2}, x_{- \beta_2} > \displaystyle{\sum_{
\theta^p \beta_1+ \beta_2=0}} \eta (p, \beta_1)
\displaystyle{\sum_{i=1}^{N}} r_i k_i (\underline{r}+ \underline{s}, \zeta_2^m) \delta (w^{-p} \zeta_1/\zeta_2).$$
For $E_4$ we use Proposition 2.3 (2) $(a=w^{-p})$ and Proposition 2.1 (2) (b).
Thus we get
$$E_4= \frac{1}{m^2} < x_{\beta_2}, x_{-\beta_2} > 
\displaystyle{\sum_{\theta^p \beta_1+ \beta_2=0}} \eta (p, \beta_1) 
k_0 (\underline{r}+ \underline{s}, \zeta_2^m) D \delta (w^{-p} 
\zeta_1/\zeta_2) .$$  
$$\frac{- 1}{m^2} < x_{\beta_2} ,x_{- \beta_2}> \displaystyle{\sum_{\theta^p
\beta_1+ \beta_2=0}} \ \eta (p, \beta_1)\frac{m}{k} \displaystyle{\sum_{i \neq 0}} (\beta_1)_i \otimes t^i (w^p \zeta_2)^i k_0 (\underline{r}+ \underline{s},
\zeta_2^m) \delta (w^{-p} \zeta_1/\zeta_2)$$
We will use the fact that $\displaystyle{\sum_{i \neq 0}} (\beta_1)_i 
\otimes t^i (w^{p}\zeta_2)^i= \beta_2 (\zeta_2) - (\beta_2)_0$
and the fact that
$$\frac{1}{k} \beta_2 (\zeta_2) k_0 (\underline{r}+ \underline{s}, \zeta_2^m)=
\beta_2 (\underline{r}+ \underline{s}, \zeta_2) \leqno{(2.6)}$$
So we get
$$E_4 = \frac{1}{m^2} <x_{\beta_2}, x_{- \beta_2}> 
\displaystyle{\sum_{\theta^p \beta_1+ \beta_2=0}}
\eta (p, \beta_1) k_0 (\underline{r}+ \underline{s}, \zeta_2^m)  
D \delta (w^{-p} \zeta_1 / \zeta_2)   $$
$$- \frac{1}{mk} <x_{\beta_2}, x_{-\beta_2} > \displaystyle{\sum_{\theta^p 
\beta_1+\beta_2 =0}} \ \eta (p, \beta_1) \displaystyle{\sum_{i \neq 0}}
(\theta^p \beta_1)_i \otimes t^i \zeta_2^i k_0  (\underline{r}+ \underline{s}, \zeta_2^m) \delta (w^{-p} \zeta_1/\zeta_2)$$
$$ =\frac{1}{m^2} <x_{\beta_2}, x_{-\beta_2} > \displaystyle{\sum_{\theta^p 
\beta_1+ \beta_2=0}} 
\eta (p, \beta_1) k_0 (\underline{r}+ \underline{s}, \zeta_2^m)
D (\delta w^{-p} \zeta_1/\zeta_2) \
$$
$$ \frac{-1}{mk} <x_{ \beta_2}, x_{-\beta_2} > \displaystyle{\sum_{\theta^p 
\beta_1+ \beta_2=0}} \ \eta (p, \beta) \theta^p \beta_1 
 (\zeta_2) k_0(\underline{r}+\underline{s},
\zeta_2^m) \delta ( w^{-p} \zeta_1/\zeta_2)$$
$$+ \frac{1}{mk} <x_{ \beta_2}, x_{-\beta_2} > \displaystyle{\sum_{\theta^p 
\beta_1+ \beta_2=0}} \ \eta (p, \beta_1) (\beta_1)_0  
k_0(\underline{r}+\underline{s},
\zeta_2^m) \delta  (w^{-p} \zeta_1/\zeta_2).
$$
Note that the second term is $-E_2$ by (2.6). Since $\theta^p \beta_1+ \beta_2=0$ we have $(\beta_1)_0+(\beta_2)_0=0$.  Thus $E_2+E_4$
$$=\frac{1}{m^2} < x_{\beta_2}, x_{-\beta_2} > \displaystyle{\sum_{\theta^p 
\beta_1+\beta_2 =0}} \ \eta (p, \beta_1) k_0 (\underline{r}+\underline{s},
\zeta_2^m) D \delta (w^{-p}  \zeta_1/\zeta_2)$$
$$- \frac{1}{mk} < x_{\beta_2}, x_{-\beta_2} > \displaystyle{\sum_{\theta^p 
\beta_1+\beta_2 =0}} \ \eta (p, \beta_1) (\beta_2)_0k_0 (\underline{r}+
\underline{s},
\zeta_2^m) \delta (w^{-p}  \zeta_1/\zeta_2).$$
Now  adding $E_1, E_2, E_3$ and $E_4$ we get the desired result.
Thus we proved that $\Omega (V) \in \underline{D}_k$.

Conversely assume that $W \in \underline{D}_k$.  Let $V= M(k) \otimes W$.
Define $X_{\alpha}( \underline{r}, \zeta) = E^- (-\alpha, \zeta) E^+
(- \alpha,\zeta) \otimes Z (\alpha, \underline{r}, \zeta)$.
$$\beta (\underline{r}, \zeta)  = \frac{1}{k} \beta (\zeta) k_0
(\underline{r}, \zeta^m).$$
The central elements to be same.  Since $ W  \in \underline{D}_k$.  The
operators $X_{\alpha} (\underline{r}, \zeta)$ and $k_0
(\underline{r}, \zeta)$ satisfy
$$ \frac{1}{k} X_{\alpha} (\underline{r}, \zeta) k_0 (\underline{s},
\zeta^m) =X_{\alpha} (\underline{r}+ \underline{s}, \zeta) \ {\rm for \
all} \ X_{\alpha} \in  {\cal G}.$$
Conditions 4 to 7 of Proposition (1.6) are easily satisfied as the
corresponding conditions are satisfied for $Z$-operators.

Condition (3) can be proved exactly as in the proof of Proposition
5.3 of [LW]. Condition (1)  is satisfied as the same relation holds
for $Z$-operators and $Z$ operator commutes with $E^{\pm}$-operators.

Consider $[\beta_1 (\zeta_1), \beta_2 (\zeta_2)]= \frac{1}{m^2} k
\displaystyle{\sum_{p \in \Z_m}} < \theta^p \beta_1, \beta_2 > D \delta
(w^{-p} \zeta_1/\zeta_2)$. 
See Theorem 2.4 of [LW].
$$
\begin{array}{lll}
[\beta_1 (\underline{r}, \zeta_1), \beta_2 (\underline{s},
\zeta_2)]&=& \frac{1}{k^2}[\beta_1 (\zeta_1), \beta_2 (\zeta_2)] k_0
(\underline{r}, \zeta^m_1) k_0 (\underline{s}, \zeta_2^m) \\
&=&\frac{1}{km^2} \displaystyle{\sum_{p \in \Z_m}} < \theta^p
\beta_1, \beta_2> k_0 (\underline{r}, \zeta_1^m) k_0 (\underline{s},
\zeta^m_2) D \delta
(w^{-p} \zeta_1/\zeta_2) 
\end{array}
$$
Now by Proposition 2.3 (2) we have
$$
= \frac{1}{m^2} \displaystyle{\sum_{p \in \Z_m}} < \theta^p \beta_1,
\beta_2> k_0 (\underline{r}+\underline{s}, 
\zeta^m_1)  D \delta (w^{-p} \zeta_1/\zeta_2) $$
$$ -\frac{1}{m^{2}k} \sum < \theta^p \beta_1, \beta_2>
D_1 k_0(\underline{r}, \zeta_1^m) \mid_{\zeta_1 = w^p \zeta_2} k_0 
(\underline{s}, \zeta_2^m) \delta( w^{-p} \zeta_1 /\zeta_2) $$
$$=\frac{1}{m^2} \displaystyle{\sum_{p \in \Z_m}} < \theta^p \beta_1,
\beta_2> k_0 (r+s, \zeta_2^m) D \delta (w^{-p} \zeta_1/\zeta_2)$$
$$+\frac{1}{m} \displaystyle{\sum_{p \in \Z_m}} < \theta^p
\beta_1, \beta_2> \displaystyle{\sum_{i=1}^{N}} r_i k_i (\underline{r}+\underline{s}, \zeta_2^m) \delta (w^{-p} \zeta_1/\zeta_2) $$
Thus we have proved the following:

\paragraph{Proposition 2.6}  The category $\underline{C}_k$ of $\tilde{L}
({\cal G}, \theta)$-modules are equivalent to the category $\underline{D}_k$ of
$Z_k$-modules. 

\section*{Section 3 (Homogeneous picture)} 

In this section our aim is to construct a faithful representation
for the untwisted toroidal Lie algebra $\tilde{\tau}$ coming from
simple, simply connected Lie-algebra ${\cal G}$.  (First note that on any
representation $\tau$ where centre $\Omega_A / d_A$ acts faithfully,
then $\tilde{\tau}$ acts faithfully). That is we are giving a
realization.  This recovers the main result of [EM].  For this we
give a representation for the $Z-$algebra such that the centre acts
faithfully. Thus we have a faithful  representation
for the toroidal Lie algebra $\tau$.

We take the automorphism $\theta = id$.  We first give a presentation for
the Lie-algebra ${\cal G}$.  Let $\stackrel{\circ}{Q}$ be the root lattice
spanned by simple roots.  The nondegenerate form is chosen  so that
$(\alpha, \alpha) =2$ for a highest root $\alpha$.  
Then it is known that $$\Phi=\{\alpha \in \stackrel{\circ}{Q} \mid (\alpha, \alpha) =2.\}$$ 
The following cocycle  on $\stackrel{\circ}{Q} \times 
\stackrel{\circ}{Q}$ is known to exists.
$$\epsilon : \stackrel{\circ}{Q} \times \stackrel{\circ}{Q} \to \{\pm
1\}$$ 
\paragraph*{3.1}
\begin{enumerate}
\item[(1)] $\epsilon (\alpha, \alpha) = (-1)^{ \frac{(\alpha, \alpha)}{2}}
$ 
\item[(2)] $ \epsilon (\alpha, \beta) \epsilon (\beta, \alpha) =
(-1)^{(\alpha, \beta)}$ 
\item[(3)] $ \epsilon (\alpha+ \beta, \gamma) = \epsilon ( \alpha,
\gamma) \epsilon (\beta, \gamma)$
\item[(4)] $\epsilon (\alpha, \beta+ \gamma) = \epsilon (\alpha, \beta)
\epsilon (\alpha, \gamma)$
\end{enumerate}

Note that $\epsilon (\alpha, \alpha)= \epsilon (\alpha, - \alpha) = -1$
for $\alpha \in \Phi$.   Then there exist vectors $x_{\alpha}, h_{\alpha}$
in ${\cal G}, \ \alpha \in \Phi$ satisfying the following:

\paragraph{(3.2)}
$$
\begin{array}{lll}
(1) ~~~ [x_{\alpha}, x_{\beta}] &=&\epsilon (\alpha, \beta)
x_{\alpha+\beta},\ {\rm if} \ \alpha+ \beta \in \Phi \\
(2) ~~~ [x_{\alpha}, x_{\beta}] &=&0 \ \ {\rm if} \ \alpha+ \beta
\notin \Phi \cup \{0 \} \\ 
&=& \epsilon (\alpha, -\alpha) h_{\alpha} \ {\rm if} \ \alpha + \beta
=0 \\
(3) ~~~ [ h_{\alpha}, h_{\beta}] &=& 0 \\
(4) ~~~ \left[ h, x_{\alpha} \right] &=& \alpha (h) x_{\alpha} ~~~
\forall h \in \underline{h}. 
\end{array}
$$

Let $\Gamma$ be a $\Z$-lattice spanned by $\alpha_1, \cdots,
\alpha_{\ell}, \delta_1, \cdots, \delta_N, d_1, \cdots, d_N$.  Define
a non-degenerate bilinear form on $\Gamma$ extending the one on
$\stackrel{\circ}{Q}$ by 
$$
\begin{array}{lll}
(\stackrel{\circ}{Q}, \delta_i) &=& (\stackrel{\circ}{Q}, d_i )=0 \\
(d_i, d_j) &=& (\delta_i, \delta_j) =0 \\
(\delta_i, d_j) = \delta_{ij}
\end{array}
$$
Any vector which is integral linear combination of $\delta_i$ is
called a null root.  For $\underline{r} \in \Z^N$,  define
$\delta_{\underline{r}} = \sum r_i \delta_i$.  Note that
$(\delta_{\underline{r}},\delta_{\underline{s}})= 0$ .  Let $Q$ be
the sub lattice spanned by $\stackrel{\circ}{Q}$ and $\delta_1,
\cdots, \delta_N$.  Extend the co-cycle $\epsilon $ to $Q$ by $\epsilon
(\alpha, \delta_{\underline{r}})=1$ for $\alpha \in Q$.  Now extend
$\epsilon$ to $Q \times \Gamma$ to be bimultiplicative in any convenient way.
Consider the group algebra $\C[\Gamma]$ and make $\C[\Gamma]$ a $\C
[Q]$ module by the following multiplication.
$$e^{\alpha} \cdot e^{\gamma} = \epsilon (\alpha, \gamma) e^{\alpha+
\gamma}.$$
Let $\overline{h}=Q  \otimes_{\Z} \C$.  Let $\overline{h}_{\pm} =
\displaystyle{\oplus_{n \stackrel{>}{<} 0}} \overline{h} \otimes t^n$.  
Consider the Fock space
$$V(\Gamma) = S (\overline{h}_-) \otimes \C [\Gamma].$$
Define operators $\alpha(0)$ on $ V (\Gamma)$ by
$$\alpha (0) \cdot u \otimes e^{\gamma} = (\alpha, \gamma) u \otimes
e^{\gamma}, \alpha \in Q.$$ 
For $\delta$ nullroot
$$\delta (n) u \otimes e^{\gamma} = \delta (n) u \otimes e^{\gamma},
n \neq 0$$
$\delta (n) u$ is multiplication if $n<0$ and differentiation if $n >0$.
This is the standard Fock space representation of $\overline{h}
\otimes \C [t, t^{-1}]$ on $V[\Gamma]$.

For a null root $\delta$ define $E^{\pm} (\delta, \zeta) = {{\rm
exp}} \displaystyle{\sum_{n \stackrel{>}{<}0}} \frac{ \delta {(\pm
n)}}{\pm n}  \zeta^{\pm n}$.  Define
operators 
$$\zeta^{\alpha (0)} u \otimes e^{\gamma} =\zeta^{(\alpha,
\gamma)} u \otimes e^{\gamma}, \alpha \in Q.$$
Consider the vertex operator
$$X (\delta, \zeta) = E^- (\delta, \zeta) \zeta^{\delta (0)} E^+ (\delta, \zeta).$$
Let $k_i (\delta, \zeta)= \delta_i (\zeta)X (\delta, z) $ for $ 1 \leq i \leq
N$ and $k_0 (\delta, z) = X (\delta, \zeta )$. From [EM] it is known that 
each $k_i (\delta, z)$ acts non trivially
and $Dk_0 (\delta_{\underline{r}}, \zeta)+\displaystyle{\sum_{i=1}^{N}} r_i k_i (\delta_{\underline{r}}, \zeta) =0$.
 Further any relation among $k_i (\delta, \zeta)$ is the one given above.
 (see Lemma $C$ of [EM]).
 
 We will now define $Z$ operators.  Define $Z (\alpha, 0, \zeta) =
\zeta^{\frac{(\alpha, \alpha)}{2}} \zeta^{-\alpha (0)} e^{\alpha}, \alpha
\in \Phi$.  Then define $ Z (\alpha, \underline{r}, \zeta) = Z
(\alpha, 0, \zeta) k_0 (\underline{r}, \zeta)$. $d_0, d_1, \cdots, d_N$
are defined naturally as grading on $Z (\alpha, \underline{r},
\zeta)$. 

We will now check the relation at (1.10) for the above $Z$-operator.
(1) to (6) are clearly satisfied from definition.  We will rewrite
the relation (7) using the fact that $\theta =Id$ and $m=1$. 
Notice also $X(0,\zeta)=1$ and hence $k_0$ acts as 1 so that $k=1.$
\paragraph*{(3.3)}
$$(1- \zeta_1/\zeta_2)^{< \beta_1, \beta_2>}Z (\beta_1, \underline{r}, \zeta_1)
Z(\beta_2, \underline{s}, \zeta_2) - (1- \zeta_2/\zeta_1)^{< \beta_1,
\beta_2>}Z (\beta_2, \underline{s}, \zeta_2) Z (\beta_1, \underline{r},
\zeta_1)$$
$=\begin{cases} \epsilon (\beta_1, \beta_2) Z (\beta_1+\beta_2,
\underline{r}+\underline{s},\zeta_2) \delta (\zeta_1/\zeta_2)\ {\rm if} \ \beta_1+\beta_2 \in \Phi\\
-<x_{\beta_2}, x_{-\beta_2}> \beta_2  k_0 (\underline{r}+ \underline{s}, \zeta_2) \delta (\zeta_1/\zeta_2)\\
+<x_{\beta_2},  x_{-\beta_2}> (\displaystyle{\sum_{i=1}^{N}} r_i k_i (\underline{r}+\underline{s}, 
\zeta_2). \delta (\zeta_1/\zeta_2)+ k_0 (\underline{r}+ \underline{s}, \zeta_2) D \delta (\zeta_1/\zeta_2)) \ {\rm if} \ \beta_1+\beta_2=0\\
=0 \ {\rm if} \ \beta_1+\beta_2 \notin \Phi \cup \{0\}.\end{cases}$
Suppose $\underline{r}=0$ and $\underline{s}=0$ then (3.3) follows from
Theorem 5.3 of $[LM]$.  For general $\underline{r}$ and $\underline{s}$,
consider left hand side of (3.3) which is equal to
$$(1- \zeta_1/\zeta_2)^{< \beta_1, \beta_2>} Z (\beta_1, 0, \zeta_1) k_0(\underline{r}, \zeta_1) 
Z (\beta_2,0, \ \zeta_2) k_0 (\underline{s}, \zeta_2) $$
$$-(1- \zeta_2/\zeta_1)^{< \beta_1, \beta_2 >} Z (\beta_2, 0, \zeta_2) k_0
(\underline{s}, \zeta_2) .Z (\beta_1, 0, \zeta_1) k_0 (\underline{r}, \zeta_1) $$

$$= k_0 (\underline{r}, \zeta_1) k_0 (\underline{s}, \zeta_2) Z (\beta_1,\beta_2) \ \mbox{where}$$ 
$Z(\beta_1,\beta_2)=\begin{cases}\epsilon (\beta_1, \beta_2) Z (\beta_1+ \beta_2, 0, \zeta_1)  
\delta (\zeta_1 / \zeta_2) \ {\rm if} \ \beta_1+ \beta_2 \in \Phi\\
 < x_{\beta_2}, x_{- \beta_2} > (- \beta_2 \delta (\zeta_1 / \zeta_2)+
D \delta  (\zeta_1/\zeta_2) \ {\rm if} \ \beta_1+ \beta_2 =0\\
0 \ {\rm if} \ \beta_1+ \beta_2 \notin  \Phi \cup \{0 \}.\end{cases}$\\
This follows from case $\underline{r} =0=s$.  Now the case $\beta_1+
\beta_2 \in \Phi $, (3.3) follows from Proposition 2.3 (1).  The case
$\beta_1+ \beta_2 \notin \Phi \cup \{0\}$ is very standard as \newline $<\beta_1, \beta_2> \geq 0.$  For the case $\beta_1+
\beta_2 =0$,$-< x_{\beta_2}, x_{-\beta_2}>$  $ k_0 (\underline{r}, \zeta_1) k_0
(\underline{r}, \zeta_2) \beta_2 \delta (\zeta_1 / \zeta_2)$ is
equal to the first term of 3.3 which follows from Proposition 2.3(1). \\
Now \\ 
$< x_{\beta_2}, x_{- \beta_2} > k_0 (\underline{r}, \zeta_1) k_0
(\underline{s}, \zeta_2).  D \delta (\zeta_1/\zeta_2)$\\
$ =< x_{\beta_2}, 
x_{- \beta_2} >( k_0 (\underline{r}+ \underline{s}, \zeta_2) D 
\delta (\zeta_1/ \zeta_2)+  \sum r_i k_i (\underline{r}+ \underline{s}, 
\zeta_2) \delta (\zeta_1 / \zeta_2)$. \\
By proposition 2.3(2). This completes the proof of (3.3).

\section*{Section 4 ~~Principal realization}

Recall that ${\cal G}$ is simple finite dimensional Lie algebra and $<,>$ a non-degenerate 
bilinear form an ${\cal G}$ . Let $\eta$ be a finite order automorphism of order p. Consider 
the affine Lie algebra ${\cal G} \otimes \C [t,t^{-1}]\oplus \C C $ with Lie bracket 
$[x \otimes t^{m}, y \otimes t^{n}]=[x,y] \otimes t^{m+n}+ \frac{<x,y>}{p} m \delta_{m+n,0} C.$
Let ${\cal G}(\eta)$ be the corresponding twisted affine Lie algebra. See $[Ka]$ for details. \\

Let $\pi$ be an automorphism of order  $K(=1,2$ or $3$) induced by an automorphism of
the Dynkin diagram of ${\cal G}$ with respect to some Cartan 
subalgebra $h$ of ${\cal G}$.  Let $\epsilon$ be $K$-th primitive root.  We
will now extend the automorphism $\pi$ to ${\cal G} \otimes A \oplus
\Omega_A / d_A= \tau$ by
\paragraph*{(4.1)}  
$$\pi (x t^{r_0} t^{\underline{r}}) = \epsilon^{-r_0} \pi (x) t^{r_0}
t^{\underline{r}} $$
$$\pi ( t^{r_0} t^{\underline{r}}{k_i}) = \epsilon^{-r_0} t^{r_0}
t^{\underline{r}}{k_i}, 
\ \  0 \leq i \leq N. $$
The aim of this section is to prove that $\overline{L} ({\cal G}, \pi)
\cong \overline{L} ({\cal G}, \theta)$ where $\theta$ is a special automorphism
depending on $\pi$.  This is a generalisation of the standard
principal realization of affine Lie-algebras given in [KKLW].

To do this we first have to define the automorphism $\theta$. For $i \in
\Z$, let ${\cal G}_{[i]}$ be the $\epsilon^i$ eigenspace of ${\cal
G}$.  Then the  fixed point space ${\cal G}_{[0]}$ is a simple
Lie-subalgebra of ${\cal G}$ and 
${\cal G}_{[0]}$ 
module ${\cal G}_{[1]}$ and ${\cal G}_{[-1]}$ are irreducible and
contragradient. 

Fix a Cartan subalgebra $\underline{t}$ of of ${\cal G}_{[0]}$ inside 
$\underline{h}$.  Let $H_j, E_j, F_j \ (1 \leq j \leq \ell)$ be a
corresponding 
set of canonical generators of ${\cal G}_{[0]}$.  Let $E_0$ be the lowest
weight vector of ${\cal G}_{[0]}$ module ${\cal G}_{[1]}$, and let $F_0$ be the
highest weight vector ${\cal G}_{[0]}$-module ${\cal G}_{[-1]}$, normalised so
that $[H_0, F_0]=2 F_0$ where $H_0 = [E_0, F_0]$.  Let $\psi_1,
\cdots, \psi_{\ell} \in \underline{t}^*$ be simple roots of
${\cal G}_{[0]}$, and let $\psi_0 \in \underline{t}^*$ be the
lowest weight 
of the ${\cal G}_{[0]}$ module ${\cal G}_{[1]}$.  For $i,j =0,1, \cdots, \ell$
set $ A_{ij} = \psi_j (H_i)$.  Then it is known that $A=(A_{ij})$ is
an indecomposable  affine Cartan matrix (see [LW] and [KKLW]).  Let $a_0,
\cdots, a_{\ell}, a_0^1, \cdots, a_{\ell}^1$  be  positive
integers such that

\paragraph*{(4.2)}
$$
\begin{array}{lll}
A(a_0, \cdots a_{\ell})^T &=& 0 \\
(a_0^1, \cdots, a_{\ell}^1) A &=& 0 \ {\rm and} \\
g.c.d (a_0, \cdots, a_{\ell}) &=& 1 \\
g.c.d (a_0^1, \cdots , a_{\ell}^1 ) &=& 1.
\end{array}
$$

Then $a_0, a_1, \cdots, a_{\ell}$ are precisely the indices of the
Dynkin diagram of $A$.  (see Table $K$ of [KKLW]).

\subsection*{(4.3)}

Note that from above tables we see that $a_0=1$ always.  Recall from [LW] that
\paragraph*{(4.4)}
$$
\begin{array}{lll}
\displaystyle{\sum_{j=0}^{\ell}} a_j \psi_j &=&0 \ {\rm and} \\
\displaystyle{\sum_{j=0}^{\ell}} a_j^1 H_j &=& 0 
\end{array}
$$
\paragraph*{(4.5) Proposition (KKLW)}  The Lie subalgebra ${\cal G} (\pi)$  of ${\cal G}
\otimes \C[t, t^{-1}] \oplus \C C$ generated by $E_i, F_i, H_i (1 \leq
i \leq \ell), E_0 \otimes t, F_0 \otimes t^{-1}, H_0$ and $C$ is
isomorphic to the affine Lie algebra corresponding to $A$.

Let $\underline{s}= (s_0, \cdots, s_{\ell})$ be a sequence of
non-negative integers, not all  $0$.
$${\rm Take} \ m= K \displaystyle{\sum_{j=0}^{\ell}} s_j a_j.$$
Define an automorphism $\theta$ of  ${\cal G}$ by the condition
$$\theta H_i = H_i, \theta E_j= w^{s_j} E_j.$$
where $w$ is m th root of unity.\\
Then $\theta$ defines an automorphism of ${\cal G}$ of order $m$.  See [LW].
Then from section 1 we can define $\overline{L} ({\cal G}, \theta)$.  Our aim
in this section is to prove $\overline{L} ({\cal G}, \pi) \cong
\overline{L} ({\cal G}, \theta)$ which is what we call principal realization of
toroidal Lie algebras.

\paragraph*{(4.6)}

Note that the  ${\cal G}$-invariant bilinear form $<,>$
on ${\cal G}$ is necessarily $\theta$ and $\pi$ invariant.  This form 
$<,>$ remains non-singular on the Cartan subalgebra $\underline{t}$ of 
${\cal G}_{[0]}$. See section 8 of [LW].

Using the restricted form we identify $\underline{t}$ and
$\underline{t}^*$.  We normalise the form $<,>$ such that

\paragraph*{(4.7)}
$$< \psi_0, \psi_0>= \frac{2 a_o^1}{K}.$$
Then we have $a_j^1=K < \psi_j, \psi_j> a_j/2$ for $j=0, \cdots,
\ell$.  (see [LW]).

Now we have the following $\underline{s}$-realization of the affine Lie
algebra ${\cal G}(\pi)$ from [KKLW].

\paragraph*{(4.8) ~~ Proposition } ([LW], [KKLW])

Let $e_j= E_j \otimes t^{s_j}, f_j = F_j \otimes t^{-s_j}, \ h_j = H_j
\otimes 1 + 2 s_j <\psi_j, \psi_j>^{-1} m^{-1} C$ inside
${\cal G} (\theta)$. Then there is an isomorphism of affine
Lie-algebras $\varphi: {\cal G} (\pi) \to {\cal G} (\theta)$ defined by
$$
\begin{array}{lll}
\varphi (E_i \otimes 1) &=& e_i ~~ 1 \leq i \leq \ell \\
\varphi (F_i \otimes 1) &=&  f_i, 1 \leq i \leq \ell \\
\varphi (E_0 \otimes t) &=& e_0 \\
\varphi (F_0 \otimes t^{-1}) &=& f_0 \\
\varphi (H_i) &=& h_i ~~ 1 \leq i \leq \ell \\
\varphi (H_0+ \frac{< E_0, F_0>}{K} C) &=& h_0.
\end{array}
$$
Further $e_i, f_i, h_i \ (0 \leq i \leq \ell)$ forms a  set of canonical
generators for the affine Lie-algebra ${\cal G} (\theta)$. Here the Lie
bracket ${\cal G} (\pi)$ is defined by the bilinear form $\frac{1}{K} < >$ 
and the Lie bracket in ${\cal G} (\theta)$ is defined by $\frac{1}{m} <,>.$

As we are interested in the principal realization we take
$\underline{s}=(1,1, \cdots, 1)$.
\paragraph*{(4.9) ~~ Remark.}  The $\pi$-invariants of
$\underline{h}$ equal to $\underline{t}$.  In particular they are
spanned by $H_i, 1 \leq i \leq \ell$.

\paragraph*{(4.10) ~~ Proposition}  Let $w$ be the Chevalley involution
automorphism of ${\cal G}$.  Let $\varphi: {\cal G} (\pi) \to {\cal G}(\theta)$ be the
isomorphism of Lie-algebras given earlier.  Then the following hold.
\begin{enumerate}
\item[(1)] $\varphi (C) =C$
\item[(2)] $\varphi (x_{\alpha} \otimes t^{r_0}) = x_{\alpha} \otimes
t^{N(\alpha) + \frac{m}{k} r_o} $ where $x_{\alpha} \otimes t^{r_0}$
is a real root vector of ${\cal G} (\pi)$.  $N(\alpha)$ is an integer independent
of $r_0$ but depends on $\alpha$.
\item[(3)] Let $ h_{\alpha}= [x_{\alpha}, w
(x_{\alpha})]$ then
$$\varphi (h_{\alpha} \otimes t^{r_0}) = h_{\alpha} \otimes
t^{\frac{m}{k} r_0} +< x_{\alpha}, w (x_{\alpha})>
\frac{N(\alpha)}{m} \delta_{r_0,0}$$
\item[(4)] $<x_{\alpha}, x_{\beta}> \neq 0$ implies $N(\alpha) +N
(\beta) =0$ 
\item[(5)] $<x_{\alpha}, x_{\beta}> = 0$ then
\end{enumerate}
$$\varphi ([x_{\alpha}, x_{\beta}] \otimes t^{r_0+s_0}) =
[x_{\alpha}, x_{\beta}] \otimes t^{(r_0+s_0) \frac{m}{k}+N
(\alpha)+N(\beta)}.$$ 
Here $[x_{\alpha}, x_{\beta}]$ could be part of real or imaginary
root. 

\paragraph*{Proof}  Let $x_{\alpha} \otimes t^{r_0}$ be a real root
vector of ${\cal G} (\theta)$.  Then ${\cal G} (\theta)$ is spanned by $x_{\alpha} \otimes^{r_0},
[x_{\alpha}, x_{\beta}] \otimes t^{r_0}$ and $C$.  We have from (4.4)
$$H_0 =- \displaystyle{\sum_{j=1}^{\ell}} \ \frac{a_j^1}{a_{o}^{1}}
H_j.$$
From Proposition (4.8) we have 
$$H_0 \otimes 1+2 <\psi_0, \psi_0>^{-1} m^{-1} C=h_0=[E_0t,
F_0t^{-1}]=[E_0,F_0]+< E_0, F_0> \frac{1}{m} C.$$
\paragraph*{(4.11)}

This implies $<E_0, F_0> = \frac{2}{< \psi_0, \psi_0>}$.  Consider
$$
\begin{array}{lll}
\varphi (H_0) &=&- \displaystyle{\sum_{j=1}^{\ell}} \frac{a_j^1}{a_0^1} \varphi (H_j) \\
&=&- \displaystyle{\sum_{j=1}^{\ell}} \frac{a_j^1}{a_0^1} ( H_j+
\frac{2}{<\psi_j, \psi_j>} \frac{1}{m} C) \\
&=& H_0 - \displaystyle{\sum_{j=1}^{\ell}} \frac{a_j^1}{a_0^1}
\frac{1}{m} \frac{a_j}{a_j^1} K C \\
&=&H_0 - \displaystyle{\sum_{j=1}^{\ell}} \frac{a_j KC}{a_0^1 m} \\
&=& H_0 - \frac{(m-a_0k)}{a_o^1m} C \\
{\rm But} \ \varphi (H_0+ \frac{<E_0, F_0>}{K}  C) &=& H_0+ \frac{2}{<\psi_0,
\psi_0>} \frac{1}{m} C
\end{array}
$$
(by Proposition 4.8)
$$
\begin{array}{lll}
\varphi (C) \frac{2}{K < \psi_0, \psi_0 >} &=& H_0+ \frac{2}{<\psi_0,
\psi_0>} \frac{1}{m} C - \varphi (H_0) \\
&=& \frac{2}{(\psi_0, \psi_0)} \frac{1}{m} C+ \frac{1}{a_o^1} C -
\frac{a_0K}{a_0^1m}C \\
\frac{\varphi (C)}{a_0^1} &=& (\frac{K}{a_0^1} \frac{1}{m}+
\frac{1}{a_o^1} - \frac{a_0K}{a_0^1m} C)
\end{array}
$$
(We are using (4.11), (4.7)).  This implies $\varphi (C) =C$ by (4.3).
\\
{\bf (2)} ~~ Clearly $\varphi (x_{\alpha} \otimes t^{r_0}) =
x_{\alpha} \otimes t^N$ for some integer.  Write $N= N(\alpha)+ r_0
\frac{m}{K}$ for some integer $N(\alpha)$ which may depend on $r_0$.
Clearly $\varphi (w (x_{\alpha}) \otimes t^{- r_0}) = w (x_{\alpha}) \otimes t^{- r_0
\frac{m}{k} - N(\alpha)}$.  Let $h_{\alpha} = [x_{\alpha}, w
(x_{\alpha})]$.  Consider the following in $\cal G(\pi)$.
$$[x_{\alpha} \otimes t^{r_0}, w(x_{\alpha}) \otimes t^{- r_0}]=
h_{\alpha}+ \frac{<x_{\alpha}, w (x_{\alpha})> r_0 C}{K}$$
Now apply $\varphi$ both sides and the bracket takes place in ${\cal G}(\theta)$
$$\varphi (h_{\alpha}+\frac{r_0}{K} <x_{\alpha}, w (x_{\alpha})>C)
 =h_{\alpha}+ \frac{<x_{\alpha}, w (x_{\alpha})>}{m}
(N(\alpha)+ \frac{m}{K} r_0) C.$$
Since $\varphi (C)=C$ we have
$$\varphi (h_{\alpha}) = h_{\alpha} + <x_{\alpha}, w(x_{\alpha})>
\frac{N(\alpha)}{m} C.$$
As $h_{\alpha}$ is independent of $r_0$ it follows that $N(\alpha)$
does not depends on $r_0$. Now consider the following in ${\cal G} (\pi)$.
$$[x_{\alpha} \otimes t^{r_o}, w (x_{\alpha}) \otimes t^{s_0}]=
h_{\alpha} \otimes t^{r_0+s_0}+ \frac{1}{K}(x_{\alpha}, w(x_{\alpha}) r_0 \delta
_{r_0+s_0,0} C.$$
As earlier apply $\varphi$ both sides
$$\varphi (h_{\alpha} \otimes t^{r_0+s_0})=h_{\alpha} \otimes
t^{(r_0+s_0) \frac{m}{k}} + <x_\alpha, w (x_{\alpha})>
\frac{N(\alpha)}{m} \delta_{r_0+s_0,0}C.$$
This proves (3).  For (4) suppose $<x_{\alpha}, x_{\beta}> \neq 0$.
Since $<,>$is $\underline{t}$-invariant, it follows that $\alpha+
\beta$ is zero root. (root with respect to $\underline{t})$.  Thus
$[x_{\alpha}, x_{\beta}]$ is a part of imaginary root and so
$[x_{\alpha}, x_{\beta}] \in \underline{h}$.  Since $\pi$ is an
automorphism we have 
$$< \pi (x_{\alpha}), \pi (x_{\beta})> = < x_{\alpha}, x_{\beta} >
\neq 0.$$
Let $e$ be the $K$ th root of unity.\\
Let $\pi (x_{\alpha}) = e^i x_{\alpha}$ and $\pi (x_{\beta}) = e^j
x_{\beta}$.  Since $< x_{\alpha}, x_{\beta} > \neq 0$ it follows that
$i+j \equiv 0 (K)$.  Thus $[x_{\alpha}, x_{\beta}]$ is $
\pi$-invariant. Consider 
$$
\begin{array}{lll}
[x_{\alpha} t^{r_0}, x_{\beta} t^{s_0}]
& =& [x_{\alpha},x_{\beta}]t^{r_0+s_0} + \frac{1}{K}<x_{\alpha}, x_{\beta} > r_0 \delta_{r_0+s_0,0} C \\
\varphi([x_{\alpha}, x_{\beta}] t^{r_0+s_0}) 
&=& [x_{\alpha},x_{\beta}] t^{(r_0+s_0) \frac{m}{K}+ N{(\alpha)}+ N{(\beta)}} \\
&+& <x_{\alpha}, x_{\beta}> N C \ {\rm for \ some } \ N.
\end{array}
$$
Now from (3) and (4.9) it follows that
$$\varphi ([x_{\alpha}, x_{\beta}] t^{r_0+s_0}) = [x_{\alpha},
x_{\beta}] t^{(r_0+s_0) \frac{m}{K}.}$$
This forces $N(\alpha) +N(\beta) =0$. \\
(5) is clear.

Now we define an isomorphism between $\overline{L} ({\cal G}, \pi)$ and
$\overline{L} ({\cal G}, \theta)$. 

\paragraph*{(4.12) ~~ Proposition}  The following map $\varphi$
define an isomorphism from $\overline{L} ({\cal G}, \pi)$ to $\overline{L}
({\cal G}, \theta)$.
$$\begin{array}{lll}
(1) ~~ \varphi (x_{\alpha} \otimes t^{r_0} t^{\underline{r}}) &=&
x_{\alpha} \otimes t^{r_0 \frac{m}{K}+ N(\alpha)} t^{\underline{r}} \\
(2) ~~~ \varphi (h_{\alpha} t^{r_0} t^{\underline{r}}) &=& h_{\alpha}
t^{r_0 \frac{m}{K}} t^{\underline{r}} \\
&+& <x_{\alpha}, w (x_{\alpha})> \frac{N(\alpha)}{m}
t^{\frac{r_0m}{K}} t^{\underline{r}} k_0 \\ 
(3) ~~~ \varphi (t^{r_0} t^{\underline{r}} k_i) &=& t^{r_0
\frac{m}{K}} t^{\underline{r}} k_i, 0 \leq i \leq N \\
(4) ~~ {\rm if} \ <x_{\alpha}, x_{\beta}> &=& 0 \\
{\varphi} ([x_{\alpha}, x_{\beta}] t^{r_0} t^{\underline{r}}
&=& [x_{\alpha}, x_{\beta}] t^{r_0 \frac{m}{K}+ N (\alpha) +
N(\beta)} t^r.
\end{array}
$$
\paragraph*{Proof}  In view of earlier Proposition the right hand side belongs to
$\overline{L} ({\cal G}, \theta)$ except possible for (3).  For (3) $t^{r_0}
t^{\underline{r}} k_i \in \overline{L} ({\cal G}, \pi)$ which means $r_0
\equiv (K)$.  Thus $r_0 \frac{m}{K} \equiv 0 (m)$ which means $t^{r_0
\frac{m}{K}} t^{\underline{r}} k_i$ belongs $ \overline{L} ({\cal G}, \theta)$.
 The fact that $\varphi$ defines an isomorphism follows by the
corresponding isomorphism of the earlier proposition.  We will verify
one bracket.  Consider
$$
\begin{array}{lll}
 {\bf 4.13} ~~ [x_{\alpha} t^{r_0}t^{\underline{r}}, x_{\beta} t^{s_0}
t^{\underline{s}}]  &=& [x_{\alpha}, x_{\beta}]t^{r_0+s_0}
t^{\underline{r}+\underline{s}} \\
&+&\frac{<x_{\alpha}, x_{\beta}>}{K} t^{r_0+s_0} t^{\underline{r}+ \underline{s}} k_0 
+<x_{\alpha},x_{\beta}> \sum r_i t^{r_0+s_0} t^{\underline{r+s}}{k_i}
\end{array}
$$
Suppose $<x_{\alpha}, x_{\beta}> =0$.  Then the $\varphi$ of both
sides are equal.  Suppose $<x_{\alpha}, x_{\beta}> \neq 0$, then by
previous proposition it follows that $[x_{\alpha}, x_{\beta}]$ is in
$\underline{h}$ and $\pi$-invariant.  Further $N (\alpha)+ N
(\beta)=0$. 
$$
\begin{array}{lll}
[\varphi (x_{\alpha} t^{r_0}t^{\underline{r}}),\varphi (x_{\beta} t^{s_0}t^{\underline{s}})]
&=& [x_{\alpha} t^{r_0 \frac{m}{K}+N(\alpha)} t^{\underline{r}}, x_{\beta} t^{s_0 \frac{m}{K}+N(\beta)}t^{\underline{s}}] \\
&=&[x_{\alpha}, x_{\beta}] t^{(r_0+s_0) \frac{m}{K}} t^{\underline{r}+\underline{s}} \\
&&+ \frac{<x_{\alpha}, x_{\beta}>}{m} (N(\alpha)+ \frac{r_0 m}{K}) t^{(r_0+s_0) \frac{m}{K}} t^{\underline{r+s}} k_0 \\
&&+ <x_{\alpha}, x_{\beta} > \sum r_i t^{(r_0+s_o) \frac{m}{K}}t^{\underline{r+s}} k_i.
\end{array}
$$
 which is exactly equal to $\varphi$ of the right hand side of (4.13).
 

\paragraph*{(4.14) ~~ Proposition} $\overline{L} ({\cal G}, \pi)$ is the
universal central extension of $L({\cal G}, \pi)$.  Follows from Remark (2.4)
of [BK].

In the next section we give a faithful realisation to $\overline{L}
({\cal G}, \theta)$ thereby giving a realization to $\overline{L}
({\cal G}, \pi)$ where the infinite dimensional centre acts faithfully.
\vskip 5mm
\noindent
{\Large{\bf Section 5}} ~~~{\bf Principal picture. }

In this section we construct level one module for the toroidal
Lie-algebra of type $A_n^K, D_n^K$ and $E_n^K$.  What we do is to
construct representation for the $Z$ algebras where the centre acts
faithfully.  That in turn constructs module for toroidal algebras of
type ADE.  This also covers the twisted case which is new result. 

Notation as in section 4.  Consider the cyclic element $E=
\displaystyle{\sum_{i=1}^{\ell}} E_i \in {\cal G}_{(1)}$.  We make the
assumption that the $\theta$-stable Cartan subalgebra $\underline{t}$ is
the centraliser of $E$.  (See [LW] for details).

Let $t_0=t \oplus \C C \oplus \C d$.  Let $L (\lambda)$ be a basic
module for ${\cal G}(\theta)$ where $\lambda \in t_0^*$.  We
renormalize the root vector $x_{\beta}, \beta \in \Phi$ such that
$[x_{\beta}, {x}_{-\beta}] =- 2/ < \beta, \beta >$ and $\eta
(p, \beta) =1$ for all $p \in \Z_m$ and for all $\beta \in
\Phi$.  As in Theorem 8.7 of [LW] choose coset  representatives
$\beta_1, \cdots, \beta_{\ell}$ for the action $\theta$ on $\Phi$ such that $(\beta_1)_0,
\cdots,  (\beta_{\ell})_0$ is a basis  for $\underline{t}$.  Let
$C_j = \lambda \Big((x_{\beta_j})_0\Big)$. 

Since $L(\lambda)$ is a basic module, we have $\dim \Omega_{L(\lambda)}=1.$ From section 8 of $[LW]$ we have 
operators $Z(\beta, \zeta)$ acting on $\Omega_{L(\lambda)}.$
\paragraph*{Propostion (5.1)}  We have the following from Section (8) of [LW].

\begin{enumerate}
\item[(1)]  dim $\Omega_{L (\lambda)} =1$
\item[(2)] $Z  (\beta_j, \zeta) = C_j $ on $\Omega_{L(\lambda)}$ 
\item[(3)] $Z (\theta^p \beta_j, \zeta) =Z (\beta_j, w^p \zeta)$ 
\item[(4)] $ \displaystyle{\prod_{p \in \Z_m}} (1-w^{-p}
\zeta_1/\zeta_2)^{< \theta^p \beta_1, \beta_2>} Z (\beta_1,\zeta_1) Z
(\beta_2, \zeta_2)$ 
\end{enumerate}
$$-\displaystyle{\prod_{p \in \Z_m}} (1-w^{-p} \zeta_2/\zeta_1)^{<
\theta^p \beta_2, \beta_1 >} Z (\beta_2, \zeta_2) Z (\beta_1, \zeta_1)$$
$$= \frac{1}{m} \displaystyle{\sum_{\theta^p \beta_1+ \beta_2 \in \Phi}}
\epsilon (\theta^p \beta_1, \beta_2) Z (\theta^p \beta_1+ \beta_2,
\zeta_2) \delta (w^{-p} \zeta_1/\zeta_2)$$
$$-2m ^{-2}  < \beta_1, \beta_1>^{-1} \displaystyle{\sum_{\theta^p
\beta_1+ \beta_2=0}} D \delta (w^{-p} \zeta_1/ \zeta_2).$$ 

\paragraph{Proof} ~ (1), (2) and (3) follows from Section 8 of [LW].  For
that just note that $\eta (p, \beta)=1$.  (4) follows from Theorem
8.7 of [LW] as $Z$ operator defined in (2) satisfy (8.21) of [LW].  We are
also using the fact that $a_0=1$.\\

Let $\Gamma$ be $\Z$-lattice spanned by $\delta_1, \cdots \delta_N, d_1,
\cdots, d_N$ with bilinear form $(\delta_i, d_j) = \delta_{ij}$ and
$(\delta_i, \ \delta_j) = (d_{i}, d_j)=0$. Let $H = \oplus \C
\delta_i$ and $H_+ = \displaystyle{\bigoplus_{n \in \Z_+,i}} \C \delta_i(n)$
Consider the symmetric algebra $S(H_+)$.

Consider the space $V(\Gamma) = S(H_+) \otimes e^\Gamma$ where
$e^\Gamma$ group algebra.

Let
$$
\begin{array}{lll}
E^+ (\delta, \zeta^m) &= & \rm{exp} \ \sum_{n > 0} \frac{\zeta(n)}{n}
\zeta^{m n}\\
E^-(\delta, \zeta^m) &=& \rm{exp} \ \sum_{n>0} \frac{\zeta(-n)}{-n}
\zeta^{-mn} 
\end{array}
$$

which act on $S(H_+)$ and on $V(\Gamma)$.

Let $\delta(0)$ act on $V(\Gamma)$ by

$$
\begin{array}{lll}
\delta(0) u \otimes e^{r} & =& (\delta, r) u \otimes e^r\\[2mm]
\zeta^{\delta(0)} u \otimes e^r &=& \zeta^{(\delta, r)} u \otimes e^r\\[2mm]
C u \otimes e^r & =& 1 u \otimes e^r
\end{array}
$$
Consider $X( \delta, \zeta^m) = E^-(\delta, \zeta^m) \zeta^{m \delta(0)}
E^+ (\delta, \zeta^m)$\\
Let $\underline{r} \in \Z^N$ and let $ \delta_{\underline{r}} = \sum \
r_i \ \delta_i$\\
Let $k_0(\underline{r}, \zeta^m) = X (\delta_{ \underline{r}}, \
\zeta^m)$\\
Let $ \delta(\zeta^m) = \sum \ \zeta^{mn}$\\
Let $k_i(\underline{r}, \zeta^m) =\delta_i(\zeta^m) X (\delta
_{\underline{r}}, \zeta^m)$.
Then by standard argument one can prove that
$$
D k_0 (\underline{r}, \zeta^m) = -m \sum_{i=1}^N r_i k_i
(\underline{r}, \zeta^m). \leqno{(5.2)}
$$
Now define $Z(\alpha, \underline{r}, \zeta) = Z(\alpha, \zeta)
k_0(\underline{r}, \zeta)$. Now we will check all the relation define
at (1.10). Remember $k=1$.
(1)  is true by definition.  (2) is true by the fact that $X(\delta_1,
\zeta^m) X (\delta_2, \zeta^m) = X (\delta_1 + \delta_2, \zeta^m)$
(3) is just (5.2) (4), (5), (6) can be easily checked. (8) is clear,
(9) is true as $\eta(p, \beta) =1$. (10) is true by definition. Thus
it remains to prove (7). To see this multiply 4 of Prop 5.1
by $k_0 (\underline{r}, \zeta_1^m) k_0 (\underline{s}, \ \zeta^m_2)$.
Consider $k_0(\underline{s}, \zeta^m_1) k_0 (\underline{s}, \zeta^m_2) D
\delta(w^{-p} \zeta_1 / \zeta_2)$ which is equal to (from
Proposition 2.3(2)) 

$$
\begin{array}{lll}
& D \delta (w^{-p} \zeta_1 / \zeta_2) k_0(\underline{s}, \zeta^m_2) k_0(\underline{s}, \zeta_2^m)\\
&-\delta(w^{-p} \zeta_1 / \zeta_2) \Big(\frac{d}{d \zeta_1} (k_0
(\underline{r}, \zeta^m_1) k_0(\underline{s}, \zeta^m_2)
\mid_{ (\zeta_1=w^{-p} \zeta_2)}\Big) \\[2mm]
&= D \delta (w^{-p} \zeta_1 / \zeta_2) k_0 (\underline{r} +
\underline{s}, \zeta^m_2)\\[2mm]
&+ \delta (w^{-p} \zeta_1 / \zeta_2) m \sum r_i k_i(\underline{r},
\zeta^m_2) k_0(\underline{s}, \zeta^m_2)\\[2mm]
&= D \delta (w^{-p} \zeta_1 / \zeta_2) k_0 (\underline{r} +
\underline{s}, \zeta^m_2)\\[2mm]
&+ m \delta (w^{-p} \zeta_1 / \zeta_2) \sum r_i k_i ( \underline{r}
+ \underline{s}, \zeta^m_2)
\end{array}
$$
Now (7) of (1.10) follows from the earlier arguments.

\pagebreak
\begin{center}
{\bf REFERENCES}
\end{center}
\vskip 6mm

\begin{enumerate}
\item[{[AABGP]}] Allison, B.N., Azam, A., Berman, S., Gao, Y. and
Pianzola, A.
Extended Affine Lie Algebras and their root systems, {\emph Mem. Amer. Math. 
Soc.}
126 (1997), No.603, 1-122.
\item[{[AG]}] Allison, B.N. and Gao, Y., The root system and the core of
an Extended Affine Lie Algebra, {\emph Selecta Math. (N.S.)},  7 (2001),
149-212.
\item[{[B1]}] Billig Yuly, Principal Vertex Operator representations
for toroidal Lie-algebras, \emph{Journal of Mathematical Physics},
39 (1998), No.7, 3844-3864
\item[{[B2]}] Billig Yuly, An extension of Kortweg-de Vries
Hierarchy arising from representations of toroidal Lie algebras,
\emph{Journal of Algebra}, 217(1999), No.1, 40-64.
\item[{[B3]}] Billig Yuly, A category of modules for the full toroidal Lie algebra, 
\emph{Int. Math. Res. Not. 2006,} Art. Id. 68395.
\item[{[BK]}] Berman Stephen and Krylyuk Yarsolav, Universal Central
extensions of twisted and untwisted Lie algebras extended over
commutative rings, Journal of Algebra, 173(1995), 302-347.
\item[{[BY]}] Berman, S. and Yuly Billig, Irreducible representations
for toroidal Lie algebras, \emph{Journal of Algebra}, 221 (1999),
188-231. 
\item[{[EM]}] Eswara Rao, S. and Moody, R.V., Vertex representations for
$N$- toroidal Lie-algebras and a generalisation of the Virasoro
algebras, \emph{Communications of Mathematical Physics} 159 (1994), 239-264.
\item[{[ISW1]}] Iohara Kenji, Saito Yoshihisa and Wakimoto Minoru,
Notes on differential equations arising from representations of
2-toroidal lie-algebras, gauge theory and integrable modules,
\emph{Progress of Theoretical Physics} Supplement No.2, 135(1999), 166-181.
\item[{[ISW2]}] Iohara Kenji, Saito Yoshihisa and Wakimoto Minoru,
Hirota Billinear forms with 2-toroidal symmetry, Physics letters A,
254(1999), No.1-2, 37-46.
\item[{[Ka]}] Kac, V.G., Infinite dimensional Lie algebras, 3rd ed. \emph{Cambridge University Press,} 1990.
\item[{[K]}] Kassel, C., K\"ahler differentials and coverings of complex simple Lie algebras extended over a commutative algebras, \emph{J Pure Appl. Algebra} 34 (1985), 265-275.
\item[{[KKLW]}] Kac, V.G., Kazhdan, D.A., Lepowsky J., and Wilson
R.L., Realization of the basic representations of the Euclidean
Lie-algebras, \emph{Advances in Mathematics}, 42(1981), No.1, 83-112.
\item[{[LW]}] Lepowsky L and Wilson R.L., The structure of standard
modules, 1: universal algebras and the Rogers-Ramanujan identies,
\emph{Inventiones Mathematicae}, 77 (1984) 199-290
\item[{[LP]}] Lepowsky J and Primc M., Standard Modules for type one
Affine Lie-algebras, {\emph Lecture notes in Mathematics}, 1052 (1984), 194-251.
\item[{[MEY]}] Moody R.V., Eswara Rao, S. and Yokonuma T, Toroidal Lie
algebras and vertex representations, \emph{Geom. ded.} 35(1990), 283-307.
\item[{[MY]}] Morita, J. and Yoshii, Y., Universal central extensions of 
Chevalley algebras over Laurent series polynomial rings and GIM Lie algebras, 
\emph{Proc. Japan Acad. Ser}, A61 (1985), 179-181.
\item[{[T]}] Tan Shaobin, Principal Construction of the toroidal Lie
algebras of type $A_1$, \emph{Math. Zeit}, 230(1999), 621-657.
\end{enumerate}

\end{document}